\documentclass{article}

\newtheorem{theorem}{Theorem}

\title{Growth of sumsets in abelian semigroups\footnote{Supported in part
by grants from the PSC--CUNY Research Award Program
and the NSA Mathematical Sciences Program.}}
\author{Melvyn B. Nathanson\\
Department of Mathematics\\
Lehman College (CUNY)\\
Bronx, New York 10468\\
Email: nathansn@alpha.lehman.cuny.edu}
\date{}

\begin{document}
\maketitle

Let $S$ be an abelian semigroup, written additively,
that contains the identity element 0.
Let $A$ be a nonempty subset of $S.$
The cardinality of $A$ is denoted $|A|.$
For any positive integer $h$, the {\em sumset} $hA$
is the set of all sums of $h$ not necessarily distinct elements of $A$.
We define $hA = \{0\}$ if $h = 0$.
Let $A_1,\ldots, A_r,$ and $B$ be nonempty subsets of $S$,
and let $h_1, \ldots, h_r$ be nonnegative integers.
We denote by
\begin{equation}          \label{semiabel:form}
B + h_1A_1 + \cdots + h_rA_r
\end{equation}
the set of all elements of $S$ that can be represented in the form
$b + u_1 + \cdots + u_r,$ where $b \in B$ and
$u_i \in h_iA_i$ for all $i = 1,\ldots, r.$
If the sets $A_1,\ldots, A_r,$ and $B$ are finite, then the
sumset~(\ref{semiabel:form}) is finite for all $h_1,\ldots, h_r.$
The {\em growth function} of this sumset is
\[
\gamma(h_1, \ldots, h_r) = |B + h_1A_1 + \cdots + h_rA_r|.
\]

For example, let $S$ be the additive semigroup
of nonnegative integers ${\mathbf N}_0$,
and let $A_1,\ldots,A_r,$ and $B$ be nonempty,
finite subsets of $\mathbf{N}_0$,
normalized so that $0 \in B \cap A_1 \cap \cdots \cap A_r$
and $\gcd(A_1 \cup \cdots \cup A_r) = 1$.
Let $b^* = \max(B)$ and $a_i^* = \max{A_i}$ for $i = 1,\ldots, r$.
Han, Kirfel, and Nathan\-son~\cite{han-kirf-nath98,nath72f}
determined the asymptotic structure of the sumset $B+h_1A_1 + \cdots + h_rA_r$.
They proved that there exist integers $c$ and $d$
and finite sets $C \subseteq [0,c-2]$ and $D \subseteq [0,d-2]$
such that
\[
B+h_1A_1 + \cdots + h_rA_r = C \cup [c, b^* + \sum_{i=1}^r a_i^* h_i -d]
\cup \left( b^* + \sum_{i=1}^r a_i^* h_i - D  \right).
\]
for $\min(h_1,\ldots, h_r)$ sufficiently large.
This implies that the growth function is eventually
a multilinear function of $h_1,\ldots, h_r,$ that is,
there exists an integer $\Delta$ such that
\[
|B+h_1A_1 + \cdots +  h_rA_r| =  a^*_1h_1 + \cdots + a^*_rh_r + b^* + 1 - \Delta
\]
for $\min(h_1,\ldots, h_r)$ sufficiently large.
The explicit determination of the sets $C$ and $D$
is a difficult unsolved problem in additive number theory.
In the case $r=1$, it is called the
{\em linear diophantine problem of Frobenius}.
For a survey of finite sumsets in additive number theory,
see Nathanson~\cite{nath96b}.

The theorem about sums of finite sets of integers
generalizes to sums in an arbitrary abelian semigroup $S.$
We shall prove that if $A_1, \ldots, A_r,$ and $B$ are finite, nonempty
subsets of $S$, then the growth function
$\gamma(h_1,\ldots, h_r)$ is eventually polynomial,
that is, there exists a polynomial $p(z_1,\ldots,z_r)$ such that
\[
\gamma(h_1, \ldots, h_r) = |B + h_1A_1+\cdots + h_rA_r| = p(h_1, \ldots, h_r)
\]
for $\min(h_1,\ldots, h_r)$ sufficiently large.
The case $r=1$ is due to Khovanskii~\cite{khov92,khov95}.
We use his method to extend the result to the case $r \geq 2.$
The idea of the proof is to show that the growth function is 
the Hilbert function of a suitably constructed module graded 
by the additive semigroup ${\mathbf N}_0^r$ of
$r$--tuples of nonnegative integers.

We need the following result about Hilbert functions.
Let $R$ be a finitely generated $\mathbf{N}_0^r$--graded
connected commutative algebra over a field $E$.
Then $R = \bigoplus_{h \in \mathbf{N}_0^r}R_h$.
Suppose that $R$ is generated by $s$ homogeneous elements $y_1,\ldots, y_s$
with $y_i \in R_{\delta_i}$, that is,
the degree of $y_i$ is $\deg y_i = \delta_i \in \mathbf{N}_0^r$.
Let $M$ be a finitely generated $\mathbf{N}_0^r$--graded
$R$--module.  For $h = (h_1,\ldots, h_r) \in \mathbf{N}_0^r$,
we define the Hilbert function
\[
H(M,h) = \dim_E\left( M_{(h_1,\ldots,h_r)} \right).
\]
For $z = (z_1,\ldots, z_r)$, we define
\[
z^h = z_1^{h_1}\cdots z_r^{h_r}.
\]
Consider the formal power series
\[
F(M,{ z}) = \sum_{h \in \mathbf{N}_0^r} H(M,h) z^h.
\]
Then there exists a vector $\beta$ with integer coordinates
and a polynomial
$P(M, z) = P(M,z_1,\ldots,z_h)$ with integer coefficients
such that
\[
F(M,{ z}) = \frac{z^{\beta}P(M, z)}{\prod_{i=1}^s(1-z^{\delta_i})}.
\]
(This is Theorem 2.3 in Stanley~\cite[p. 33]{stan83}).

\begin{theorem}
Let $A_1,\ldots,A_r,$ and $B$ be finite, nonempty subsets of an abelian
semigroup $S$.  There exists a polynomial $p(z_1,\ldots,z_r)$
such that
\[
|B + h_1A_1 + \cdots + h_rA_r| = p(h_1,\ldots,h_r)
\]
for all sufficiently large integers $h_1,\ldots, h_r.$
\end{theorem}

{\bf Proof}.
For $i = 1,\ldots, r,$ let
\[
A_i = \left\{ a_{i,1}, \ldots, a_{i,k_i} \right\},
\]
where
\[
|A_i| = k_i \geq 1.
\]
We introduce a variable $x_{i,j}$
for each $i=1,\ldots,r$ and $j = 1,\ldots, k_i$.
Fix a field $E$.
We begin with the polynomial ring
\[
R = E[x_{1,1},\ldots,x_{r,k_r}]
\]
in the $s = k_1 + \cdots + k_r$ variables $x_{i,j}$.
The algebra $R$ is connected since it is an integral domain
(cf. Hartshorne~\cite[Exercise 2.19, p. 82]{hart77}).
For each $r$--tuple $(h_1,\ldots,h_r) \in {\mathbf N}_0^r$ we let
\[
R_{(h_1,\ldots,h_r)}
\]
be the vector subspace of $R$ consisting of all polynomials
that are homogeneous of degree $h_i$
in the variables $x_{i,1},\ldots, x_{i,k_i}$.
In particular, $E = R_{(0,\ldots,0)}$.
Then
\[
R = \bigoplus_{(h_1,\ldots,h_r)\in {\mathbf N}_0^r} R_{(h_1,\ldots,h_r)}.
\]
The multiplication in the algebra $R$ is consistent
with this direct sum decomposition in the sense that
\[
R_{(h_1,\ldots,h_r)}R_{(h'_1,\ldots,h'_r)}
\subseteq R_{(h_1+h'_1,\ldots,h_r+h'_r)},
\]
and so $R$ is graded by the semigroup ${\mathbf N}_0^r$.

Next we construct an ${\mathbf N}_0^r$--graded $R$--module $M.$
To each $r$--tuple $(h_1,\ldots,h_r) \in {\mathbf N}_0^r$ we associate
a finite-dimensional vector space $M_{(h_1,\ldots,h_r)}$ over the field $E$
in the following way.  To each element
\[
u \in B + h_1A_1 + \cdots + h_rA_r
\]
we assign the symbol
\[
[u,h_1,\ldots,h_r].
\]
Let $M_{(h_1,\ldots,h_r)}$ be the vector space
consisting of all $E$--linear combinations of these symbols.
Then
\begin{equation}    \label{semiabel:growth}
\dim_E M_{(h_1,\ldots,h_r)} = |B + h_1A_1 + \cdots + h_rA_r|.
\end{equation}
Let
\[
M = \bigoplus_{(h_1,\ldots,h_r)\in {\mathbf N}_0^r} M_{(h_1,\ldots,h_r)}.
\]
This is an ${\mathbf N}_0^r$--graded vector space over $E.$

To make $M$ a module over the algebra $R,$ we must construct
a bilinear multiplication $R\times M \rightarrow M.$
We define the product of the variable $x_{i,j} \in R$ and
the basis element $[u,h_1,\ldots,h_r] \in M$ as follows:
\[
x_{i,j}[u,h_1,\ldots,h_r]
=[u+a_{i,j},h_1,\ldots,h_{i-1},h_i+1,h_{i+1},\ldots,h_r].
\]
This makes sense since
\[
u \in B + h_1A_1 + \cdots + h_iA_i + \cdots + h_rA_r
\]
and so
\[
u + a_{i,j} \in B + h_1A_1 + \cdots + (h_i+1)A_i + \cdots + h_rA_r.
\]
This induces a well-defined multiplication of elements of $M$
by polynomials in $R$ since, if $i < i'$,
\begin{eqnarray*}
\lefteqn{ x_{i',j'}\left(x_{i,j}[u,h_1,\ldots,h_r]\right) } \\
& = & x_{i',j'}[u+a_{i,j},h_1,\ldots,h_i+1,\ldots,h_r]  \\
& = & [u+a_{i,j}+a_{i',j'},h_1,\ldots,h_i+1,\ldots,h_{i'}+1,\ldots,h_r]   \\
& = & [u+a_{i',j'}+a_{i,j},h_1,\ldots,h_i+1,\ldots,h_{i'}+1,\ldots,h_r]   \\
& = & x_{i,j}[u+a_{i',j'},h_1,\ldots,h_{i'}+1,\ldots,h_r]  \\
& = & x_{i,j}\left(x_{i',j'}[u,h_1,\ldots,h_r]  \right).
\end{eqnarray*}
The case $i \geq i'$ is similar.  Note that this is the only place
where we use the commutativity of the semigroup $S.$
It follows that $M$ is an $R$--module.  Moreover,
\[
R_{(h_1,\ldots,h_r)}M_{(h'_1,\ldots,h'_r)}
\subseteq M_{(h_1+h'_1,\ldots,h_r+h'_r)},
\]
and so $M$ is a graded $R$--module.
Furthermore, the finite set
\[
\{[b,0,\ldots,0] : b \in B\} \subseteq M
\]
generates $M$ as an $R$--module.

Since $x_{i,j} \in R_{\delta_{i,j}}$, where
$\deg(x_{i,j}) = \delta_{i,j}$ is the $r$--tuple
whose $i$--th coordinate is 1 and whose other coordinates are 0,
and since
\[
\frac{1}{(1-z_i)^{k_i}}
= \sum_{h_i=0}^{\infty} {h_i+k_i-1 \choose k_i - 1}z_i^{h_i},
\]
we have
\begin{eqnarray*}
F(M,z)
& = & \sum_{h \in \mathbf{N}_0^r} H(M,h) z^h  \\
& = & \frac{z^{\beta}P(M, z)}
{\prod_{i=1}^r\prod_{j=1}^{k_i}(1-z^{\delta_{i,j}})}  \\
& = & \frac{z^{\beta}P(M, z)}{\prod_{i=1}^r(1-z_i)^{k_i}}  \\
& = & z^{\beta}P(M,z)\prod_{i=1}^r
\sum_{h_i=0}^{\infty} {h_i+k_i-1 \choose k_i - 1}z_i^{h_1}  \\
& = & z^{\beta}P(M,z)
\sum_{h = (h_1,\ldots,h_r)\in \mathbf{N}_0^r}
\prod_{i=1}^r {h_i+k_i-1 \choose k_i - 1} z^h.
\end{eqnarray*}
This implies that the Hilbert function $H(M,h)$ is a polynomial
in $h_1,\ldots, h_r$ for $\min(h_1,\ldots, h_r)$ sufficiently large.
By~(\ref{semiabel:growth}), the growth function is the Hilbert function of $M$.
This completes the proof.

\end{document}